\documentstyle[11pt]{article}

\catcode`\@=11

\@addtoreset{equation}{section}

\catcode`\@=12

\def\aphi{{|\Phi|}}
\def\aaphi{{||\Phi||}}

\def\anu{{|\nu|}}

\def\CC{{\cal C}}
\def\FF{{\cal F}}
\def\GG{{\cal G}}
\def\HH{{\cal H}}

\def\LL{{\cal L}}
\def\MM{{\cal M}}

\def\PP{{\cal P}}

\def\SS{{\cal S}}
\def\TT{{\cal T}}

\def\MF{{\cal MF}}
\def\PMF{{\cal PMF}}

\def\PML{{\cal PML}}

\def\half{{\textstyle{1\over2}}}

\def\ep{\epsilon}

\def\union{\cup}
\def\intersect{\cap}

\def\Empty{}
\catcode`\@=11

\def\eqalign#1{\,\vcenter{\openup\jot\m@th
  \ialign{\strut\hfil$\displaystyle{##}$&$\displaystyle{{}##}$\hfil
	\crcr#1\crcr}}\,}

\long\def\@savemarbox#1#2{\global\setbox#1\vtop{\hsize\marginparwidth 
   \@parboxrestore\tiny\raggedright #2}}
\def\mydesc{\list{}{\labelwidth\z@ \itemindent-\leftmargin
\listparindent 1.5em
\let\makelabel\descriptionlabel}}

\def\Smatrix#1{\,\vcenter{\normalbaselines\baselineskip=2pt\m@th
    \ialign{\hfil$\scriptstyle ##$\hfil&&\hfil$\scriptstyle ##$\hfil\crcr
      \crcr
      #1\crcr\crcr}}\,}

\def\fnum@figure{{\small Figure \thefigure}}
\def\fakefigure{\def\@captype{figure}}

\long\def\@makecaption#1#2{
    \vskip 10pt 
    \def\FCap{#2} \def\NoCap{\ignorespaces}
    \ifx \FCap\NoCap
       \setbox\@tempboxa\hbox{#1}  % This is to avoid the damn colon.
      \else
       \setbox\@tempboxa\hbox{#1: \small \it #2}
    \fi
    \ifdim \wd\@tempboxa >\hsize   % IF longer than one line:
        \unhbox\@tempboxa\par      %   THEN set as ordinary paragraph.
      \else                        %   ELSE  center.
        \hbox to\hsize{\hfil\box\@tempboxa\hfil}  
    \fi}

\def\@oddhead{\hbox{}\rightmark \hfil \rm\thepage}% Right heading.
\def\sectionmark#1{\markright {\sc{\ifnum \c@secnumdepth >\z@
      \S\thesection.\hskip 1em\relax \fi #1}}}

\catcode`\@=12

\def\oplabel#1{
  \def\OpArg{#1} \ifx \OpArg\Empty {} \else
  	\label{#1}
  \fi}
\newtheorem{theoremSt}{Theorem}[section]

\newtheorem{exerciseSt}[theoremSt]{Exercise}

\def\MakeStEnv#1{
  \newenvironment{#1}[2]{
  \begin{#1St} \oplabel{##1}%
  \global\def\CrntSt{\thetheoremSt}%
  \def\Labl{##1}\ifx\Labl\Empty{} \else {\rm (\Labl)\,}\fi%
  {\rm ##2}%
}{ 
  \end{#1St} }
}
\MakeStEnv{theorem}
\MakeStEnv{corollary}
\MakeStEnv{proposition}
\MakeStEnv{lemma}
\MakeStEnv{define}
\MakeStEnv{example}
\MakeStEnv{conjecture}

\def\lref#1{\ref{#1} (#1)}% Prints label as well as number
\newlength{\saveu}
\def\blankspace#1#2{
  \setlength{\saveu}{#2in}
  \divide\saveu by 2
  \vskip\saveu
  \centerline {\small \mbox{}#1}
  \vskip\saveu
}

\long\def\placefig#1#2#3#4{
  \begin{figure}[#4]
    \blankspace{#1}{#2}
    \caption[#1]{#3}
    \oplabel{#1}
  \end{figure}}

\newenvironment{proof}[1]{
  \def\PfArg{#1}
  \ifx\PfArg\Empty
  	\edef\PfArg{\CrntSt}  \fi
 \startproof{\PfArg}%
}{ 
  \finishproof{\PfArg}
}
\newcommand{\startproof}[1]{
  \medbreak\mbox{}
  {\it Proof of #1:}%
}
\newcommand{\finishproof}[1]{ 
  \def\FPArg{#1}
  \ifx\FPArg\Empty
  	\def\FPArg{\CrntSt}  \fi
  \smallbreak\noindent\makebox[\textwidth]{\hfill\fbox{\FPArg}}
  \medbreak\noindent
}

\newcommand{\numreg}[3]{}		%Does nothing
\newcommand{\reals}{{\bf R}}

\def\R{\reals}

\def\bigrightarrow#1{\hbox to #1{\rightarrowfill}}
\def\bigleftarrow#1{\hbox to #1{\leftarrowfill}}

\def\implies{\Rightarrow}

\def\union{\cup}
\def\intersect{\cap}
\def\semidir{\mathrel{\hbox{\vrule depth-.03ex height1.1ex\kern-0.15em$\times$}}}

\def\nbd{{\cal N}}

\def\til{\widetilde}

\begin{document}

\author{Yair N. Minsky \\Technion IIT, Haifa , Israel 32000\\ Yale
University, New Haven CT 06520}
\title{Quasi-projections in Teichm\"uller space}

\maketitle

\def\Lbar{\overline L}
\def\Str{{\rm Stretch}}
\def\Mm{{\rm Maxmin}}
\def\mM{{\rm Minmax}}
\newlabel{Slowmaps intersection lemma}{{3.2}{9}}

Many parallels have been drawn between the geometric properties of the
Teichm\"uller space of a Riemann surface, and those of complete,
negatively curved spaces (see for example \cite{bers,kerckhoff2,wolpert}).
This paper investigates one such parallel --
the contracting properties of certain projections to geodesics.

We will use the Teichm\"uller metric throughout the paper (the other
famous metric on Teichm\"uller space, the Weil-Petersson metric, is
negatively curved although it is not complete. The Teichm\"uller
metric is complete, but not negatively curved -- see \cite{masur}.)

Every closed subset $C$ of a complete metric space $X$ determines
a ``closest-points'' projection,
defined as a map $\pi_C:X \to \PP(C)$, where $\PP(C)$ is the set
of closed subsets of $C$. Namely,
$$\pi_C(x) = \{y\in C:d(x,y) = d(x,C)\}$$
where $d(x,C) = \inf_{y\in C} d(x,y)$. Suppose now that $X$ is a
simply-connected Riemannian manifold with non-positive sectional
curvatures, and that $C$ is a geodesic segment, ray or line. In this
case $\pi_C(x)$ always contains a single point. If the sectional
curvatures are bounded above by a negative constant, then
$|d\pi_C|\to 0$ as $d(x,C)\to \infty$. In other words, $\pi$ is
contracting at large distances. 

A standard consequence of this is that ``quasi-geodesics'' are always
in bounded neighborhoods of actual geodesics in such a space (see
section \ref{applications}).

\medskip

If $X=\TT(S)$, the Teichm\"uller space of a surface $S$ of finite
type, we must be satisfied with coarser results. We show that, for a
certain class of geodesics, the images of the projection do not
deviate too much from being single points (see section
\ref{preliminaries} for definitions): 

\medskip\noindent
{\bf Theorem \ref{Quasi-map}} (Quasi-map)
{\em
Given $\ep$ there is a constant $b_0(\ep,\chi(S))$ such that for any
$\ep$-thick geodesic $L$ the projection $\pi_L$ satisfies:
$${\rm diam}(\pi_L(\sigma))\le b_0$$
for all $\sigma\in\TT(S)$. 
}
\medskip

This theorem is a consequence of a characterization of $\pi_L$ in
terms of an extremal-length ratio problem. Minimizing the distance
of $\sigma\in\TT(S)$ to a geodesic $L$ is equivalent to solving the
problem of finding
$$\inf_{\tau\in L} \sup_\lambda {E_\tau(\lambda)\over
E_\sigma(\lambda)}$$
where the supremum is over $\lambda$ in $\PMF(S)$, the set of
the set of projectivized measured
foliations of compact support in $S$. This is because the supremum of
extremal length ratios is just $\exp 2d(\sigma,\tau)$.
Our second theorem, which we state informally, is: 

\medskip\noindent
{\bf Theorem \ref{Characterize projection}} (Characterize projection)
{\em
If $\sigma\in\TT(S)$ and $L$ is $\ep$-thick,
the solutions $\tau$ of the inverted problem
$$\sup_\lambda \inf_{\tau\in L} {E_\tau(\lambda)\over
E_\sigma(\lambda)}$$
are a bounded distance apart from $\pi_L(\sigma)$. 
}

\medskip

This turns out to be of great assistance in controlling the
projection. In particular we obtain the following average contraction
theorem: 

\medskip\noindent
{\bf Theorem \ref{Contraction at a distance}} (Contraction at a
distance)
{\em
There exist constants $b_1,b_2$ depending on $\ep,\chi(S)$ such
that if $L$ is an $\ep$-thick geodesic and
$d(\sigma,L)> b_1$ then, setting $r=d(\sigma,L)-b_1$,
$${\rm diam}(\pi_L(\nbd_r(\sigma)))\le b_2.$$
}

\medskip

Here $\nbd_r$ denotes an $r$-neighborhood in the Teichm\"uller metric,
and we have taken the liberty of denoting by $\pi_L(A)$ the union
$\union_{a\in A}\pi_L(a)$.

\medskip

In section \ref{applications} we give some fairly easy consequences of
these results, which are directly analogous to well-known properties
of hyperbolic space.
The first, which was the original motivation for the
results, is the aforementioned ``stability'' of certain
quasi-geodesics: a quasi-geodesic whose endpoints are connected by an
$\ep$-thick geodesic must remain in a  bounded neighborhood of the geodesic.
The second concerns the projected images of regions in Teichm\"uller
space which are analogous to horoballs in hyperbolic space.

%%% Check this before putting it in:

%\medskip

%We should comment briefly on the restriction to ``thick'' geodesics,
%which is used very heavily in the paper. In the general case, a
%geodesic may travel through parts of $\TT(S)$ where some
%non-peripheral curves on the surface develop arbitrarily small
%extremal length, and here we do not expect the contraction properties
%to continue to hold. Even the bound on the  diameter of the projection
%is a problem, as one may see by considering a limiting case: 
%given a separating simple closed curve $\alpha$, the geometry of the
%codimension 2 slice $\{\sigma:E_\sigma(\alpha) = \ep\}$ is
%approximated, as $\ep\to0$, by the Teichm\"uller space of the
%disconnected surface obtained by deleting $\alpha$ from $S$. This
%space is a product of Teichm\"uller spaces, and the Teichm\"uller
%metric on it is a sup metric. The reader can easily check that
%closest-point projections are badly behaved in such a space. 

\section{Preliminaries}
\label{preliminaries}
We begin with a brief summary of results and notation. Let $S$ be a
surface of finite genus with finitely many punctures, and let $\TT(S)$
be the Teichm\"uller space of conformal structures of finite type on
$S$, where two structures are considered equivalent if there is a
conformal isomorphism of one to the other, which is isotopic to the
identity on $S$. By ``finite type'' we mean that neighborhoods of the
punctures are conformally equivalent to the plane minus a point,
rather than a disk (hyperbolically, these are cusps rather than
funnels). 

$\TT(S)$ has a natural topology, and is homeomorphic to a
finite-dimensional Euclidean space. The {\em Teichm\"uller distance}
between two points $\sigma,\tau\in\TT(S)$ is defined as 
$$d(\sigma,\tau) = \half\log K(\sigma,\tau),$$
where $K(\sigma,\tau)=K(\tau,\sigma)$ is the smallest possible quasi-conformal
dilatation of a homeomorphism from $(S,\sigma)$ to $(S,\tau)$ isotopic
to the identity. 

Let $\alpha$ be a homotopy class of simple closed curves in $S$.
Given $\sigma\in\TT(S)$, one may define the {\em extremal length} of
$\alpha$ in $\sigma$, as the inverse of the conformal modulus of the thickest
embedded annulus homotopic to $\alpha$ in $S$:
$$E_\sigma(\alpha) = \inf_{A\sim\alpha} {1\over{\rm Mod(A)}}.$$
This number is positive provided $\alpha$ is non-peripheral -- that
is, not deformable into a puncture. We will make heavy use of the fact
(see \cite{kerckhoff}) that
\begin{equation}
K(\sigma,\tau) = \sup_\alpha {E_\tau(\alpha)\over E_\sigma(\alpha)}.
\label{K ratio}
\end{equation}

In order to obtain this supremum as a maximum, we may enlarge
the space of homotopy classes of non-peripheral simple closed curves
to obtain the space of ``measured foliations with compact
support'' (see \cite{travaux,levitt}), $\MF(S)$. Elements of $\MF(S)$
are equivalence classes of foliations of $S$ with saddle singularities
(negative index) in the interior, and index 1/2 singularities at the
poles, equipped with transverse measures. The equivalence is via
isotopy, and saddle-collapsing (Whitehead) moves. $\MF(S)$ carries a
natural topology, and is homeomorphic to a Euclidean space of the same
dimension as $\TT(S)$. Multiplication of the measure by positive
scalars gives a ray structure to $\MF(S)$, and the quotient of
$\MF(S)$ (minus the empty foliation) by this multiplication is the
projectivized space $\PMF(S)$, which is a sphere. 

A homotopy class of simple closed curves, equipped with a real
number giving the 
measure, is represented in $\MF(S)$ as a foliation whose non-singular
leaves are in the homotopy class and fit together in a cylinder whose
height is the measure. These special foliations are dense in $\MF(S)$,
and the extremal length function extends \cite{kerckhoff} to a function
$E_\sigma(\lambda)$ which is continuous in both $\sigma$ and
$\lambda$, and scales quadratically:
$$E_\sigma(a\lambda) = a^2 E_\sigma(\lambda).$$
Moreover, the supremum in (\ref{K ratio}) is realized by a unique
projective class in $\PMF(S)$. 

Fixing $\sigma\in\TT(S)$, let $QD(S,\sigma)$ denote the space of
integrable holomorphic quadratic differentials on $(S,\sigma)$ (see
\cite{gardiner}). This
space may be identified with $\MF(S)$ by associating to a quadratic
differential $\Phi$ its horizontal foliation, which we call $\Phi_h$.
This correspondence is a homeomorphism \cite{hubbard-masur}. We also
have the equality \cite{kerckhoff}
$$E_\sigma(\Phi_h) = ||\Phi||.$$
We note also that the vertical foliation $\Phi_v$, which is the
orthogonal foliation to $\Phi_h$, is just $(-\Phi)_h$. 

Geodesics in $\TT(S)$ are determined by holomorphic quadratic
differentials (and hence measured foliations) as
follows: Fix $\sigma\in\TT(S)$ and $\Phi\in QD(S,\sigma)$. The family
$\{\sigma_t\}_{t\ge 0}$
of conformal structures obtained by contracting along the leaves of
$\Phi_h$ by $e^t$ and expanding along $\Phi_v$ by $e^{-t}$, forms a
geodesic ray in the Teichm\"uller metric. Geodesics exist and are
unique between any two points in $\TT(S)$. We also have
$$d(\sigma,\sigma_t) = t$$
and
$$E_{\sigma_t}(\Phi_h) = e^{-2t} E_{\sigma}(\Phi_h),$$
$$E_{\sigma_t}(\Phi_v) = e^{2t} E_{\sigma}(\Phi_v).$$

\medskip

Our last group of definitions has to do with intersection number. The
standard  geometric (unoriented) intersection number between two simple closed
homotopy classes may be extended to a continuous function on
$\MF(S)\times\MF(S)$, homogenous in both arguments. One may best
visualize this in the case when one of the arguments is a simple
closed curve and the other is any measured foliation. The intersection
number is the smallest possible transverse measure deposited on a
homotopic representative of the simple closed curve.

Intersection numbers may be related to extremal lengths by the
following easily proven lemma which appears for example in \cite{minsky-slow}:
\begin{lemma}{Extremal length and intersection number}{}
If $\alpha, \beta$ are any two measured foliations on any Riemann
surface $X$, then
$$E_X(\alpha)E_X(\beta) \ge i(\alpha,\beta)^2.$$
\end{lemma}

We say a foliation $\lambda$ is {\em complete} if for any other
foliation $\mu$, $i(\lambda,\mu) = 0$ implies that the underlying leaf
structures of $\lambda$ and $\mu$ are the same (up to isotopy and
Whitehead moves). 

Let $L$ denote a geodesic segment, ray or line in $\TT(S)$. We say
that $L$ is {\em $\ep$-thick}, for $\ep>0$, if for all $\tau\in L$ and
all non-peripheral simple closed curves $\gamma$ in $S$,
$E_\tau(\gamma)\ge\ep$. 

We will have need of another simple observation, which appears in \cite{masur}
and again in \cite{minsky-2d}:

\begin{lemma}{Thick implies complete}
If $\Phi_h\in\MF(S)$ determines an $\ep$-thick Teichm\"uller ray, then
$\Phi_h$ is complete. 
\end{lemma}

\section{Projecting onto thick geodesics}
\label{projecting}
Our first result is that the projection to a thick geodesic is a
``quasi-map'', in the following sense:

\begin{theorem}{Quasi-map}{}
Given $\ep$ there is a constant $b_0(\ep,\chi(S))$ such that for any
$\ep$-thick geodesic $L$ the projection $\pi_L$ satisfies:
$${\rm diam}(\pi_L(\sigma))\le b_0$$
for all $\sigma\in\TT(S)$. 
\end{theorem}

\begin{proof}{}
Let $\Lbar$ denote the complete geodesic containing $L$. Choose an
arbitrary origin and orientation on $\Lbar$ and let the parameter $t$
denote twice the signed distance from the origin (the factor of 2
is for later convenience). Then $L$ corresponds to an
interval $[a,b]\intersect \R$, where $-\infty\le a<b\le \infty$.
(We will henceforth abbreviate this to just $[a,b]$).
Denote by
$\Lbar(t)$ the conformal structure on $S$ given by the point
parametrized by $t$; or by $L(t)$ if $t\in[a,b]$.

\medskip

For each $\alpha\in\MF(S)$ and $t\in \R$, let $E_t(\alpha)$ be short for
the extremal length $E_{\Lbar(t)}(\alpha)$. 
Fixing a point
$\sigma\in\TT(S)$, define also 
$$R_{\sigma,t}(\alpha) = {E_t(\alpha)\over E_\sigma(\alpha)},$$
which we will tend to abbreviate as $R_t(\alpha)$.
Note that $R_t$ is invariant under scaling, making it a function
defined on projective classes in $\PMF(S)$. 

For each $t$ we have $d(\sigma,L(t)) = \half\log \sup_\alpha R_t(\alpha)$,
so finding $\pi_L(\sigma)$ reduces to minimizing this quantity over
$[a,b]$. 

Accordingly, define
\begin{eqnarray*}
\mM(\sigma,L) = \{(\alpha,t)\in\MF(S)\times[a,b]:
		R_t(\alpha) & = & \inf_s\sup_\beta R_s(\beta) \\
			 & = & \sup_\beta R_t(\beta)\}
\end{eqnarray*}
(where infima and suprema are over $[a,b]$ and $\MF(S)$,
respectively). Let $T_1(\sigma,L)\subset[a,b]$ denote the set of values of $t$
which appear in $\mM(\sigma,L)$, and let $A_1(\sigma,L)\subset \MF(S)$
denote the set of values of $\alpha$ which appear. $\pi_L(\sigma)$ is
exactly the set $\{L(t): t\in T_1(\sigma,L)\}$.

Now consider the dual problem,
\begin{eqnarray*}
\Mm(\sigma,L) = \{(\alpha,t)\in\MF(S)\times[a,b]:
		R_t(\alpha) & = & \sup_\beta\inf_s R_s(\beta) \\
			 & = & \inf_s R_s(\alpha)\}
\end{eqnarray*}
and similarly define $T_2$ and $A_2$. 
The main idea of this paper is to replace the first problem with the
second. In particular theorem \ref{Quasi-map} will follow from

\begin{theorem}{Characterize projection}{}
Let $L$ be an $\ep$-thick geodesic in $\TT(S)$. Then for any
$\sigma\in \TT(S)$, 
$T_1(\sigma,L)$ and $T_2(\sigma_L)$ have bounded diameters, and are a
bounded distance apart.
These bounds depend only  on $\ep$ and $\chi(S)$. 
\end{theorem}

\begin{proof}{}
The first step is to show that $E_t(\alpha)$, as a function of $t$, is
coarsely approximated by a strictly convex function. This in
particular will show that the minima of $E_t(\alpha)$ occur on
sets of bounded diameter. 

Let $\Phi$ be the holomorphic quadratic differential on $\Lbar(0)$
that determines the geodesic. The extremal lengths of the horizontal
and vertical foliations of $\Phi$ along the geodesic are
$$E_t(\Phi_h) = e^{-t}E_0(\Phi_h)$$
and
$$E_t(\Phi_v) = e^t E_0(\Phi_v).$$
We also have $E_0(\Phi_h) = E_0(\Phi_v) = \aaphi$, and we may
normalize so that $\aaphi = 1$. 

We will need to know (see \cite{gardiner-masur}) that
$\Phi_h$ and $\Phi_v$ {\em bind}, in the sense that every non-zero element of
$\MF(S)$ has positive intersection with at least one of them. 
Thus we may define the strictly positive function
$$e_t(\alpha) = \half\left({i(\alpha,\Phi_h)^2\over E_t(\Phi_h)}  + 
			   {i(\alpha,\Phi_v)^2\over E_t(\Phi_v)} \right).$$
We then have
\begin{lemma}{Extremal length estimate}{}
There is a fixed $c_0 = c_0(\ep,\chi(S))$ such that if $L$ is $\ep$-thick
then for any $\alpha\in \MF(S)$ and $t\in[a,b]$,
$$e_t(\alpha)\le E_t(\alpha) \le c_0e_t(\alpha).$$
\end{lemma}

\begin{proof}{}
The left side actually does not depend on the $\ep$-thickness -- it
follows directly from lemma \ref{Extremal length and
intersection number}, applied to $\alpha$ and each of $\Phi_h$ and
$\Phi_v$. 

The right side follows from a compactness argument.  Let $\Phi^t$
denote the image of $\Phi$ in $L(t)$ under the Teichmuller map from
$L(0)$, scaled so that $||\Phi^t||=1$. Then as measured foliations we
have $\Phi^t_h = e^{-t/2}\Phi_h$ and $\Phi^t_v = e^{t/2}\Phi_v$. Thus,
$$e_t(\alpha) = \half\left(i(\alpha,\Phi^t_h)^2 +
i(\alpha,\Phi^t_v)^2\right).$$ 
The ratio $E_t(\alpha)/e_t(\alpha)$ is therefore a continuous positive function
of the projectivized measured foliation $[\alpha]\in\PMF(S)$ and the
pair $(L(t),\Phi^t)$ in the total space of Riemann surfaces equipped with
holomorphic quadratic differentials.
Fixing the Riemann surface $L(t)$, the space $\PMF(S)\times
PQD(L(t))$ is compact; and the function is invariant under the action
of the mapping class group. Further, the image of $L(t)$ in the
moduli space is restricted to a compact set by the assumption of
$\ep$-thickness. It follows that $E_t(\alpha)/e_t(\alpha)$ is bounded.
\end{proof}{}

Fix a foliation $\alpha\in\MF(S)$, and
let $s_\alpha$ denote the point in $[a,b]$ where $e_t(\alpha)$ takes
its minimum. Since $e_t(\alpha)$ is a strictly positive sum of
exponentials, a brief calculation shows that, excluding the case
$s_\alpha=\pm \infty$, 
\begin{equation}
\half e_{s_\alpha}(\alpha) \exp|t-s_\alpha| \le e_t(\alpha) \le
2e_{s_\alpha}(\alpha) \exp|t-s_\alpha|
\label{e is exp}
\end{equation}
for $t\in[a,b]$.

Let $m(\alpha)$ denote the set of minima of $E_t(\alpha)$.
It follows immediately from (\ref{e is exp}) and
lemma \ref{Extremal length estimate}
that $m(\alpha)$ is constrained by
\begin{equation}
m(\alpha) \subseteq [s_\alpha-2\log2c_0,s_\alpha+2\log2c_0].
\label{bound malpha}
\end{equation}
Since $T_2(\sigma,L) = \bigcup_{\alpha\in A_2(\sigma,L)} m(\alpha)$,
our next task will be to control the points $s_\alpha$.

\medskip

Define
$$I_t(\alpha,\beta) = {i(\alpha,\beta)^2\over E_t(\alpha)E_t(\beta)}.$$

\begin{lemma}{Distance implies intersection}{}
There are constants $D$ and $c_1$, depending only on $\ep$ and
$\chi(S)$, such that for any $\alpha,\beta\in\MF(S)$,
$$|s_\alpha - s_\beta| \ge D \quad\implies\quad 
I_{s_\alpha}(\alpha,\beta) \ge c_1.$$
\end{lemma}

\paragraph{Remark.}
This lemma is similar to (and was motivated by) lemma \ref{Slowmaps
intersection lemma} from \cite{minsky-slow}, in which the role of $L$
is played by a hyperbolic 3-manifold $N$ homeomorphic to $S\times \R$,
and the role of the minimum point $s_\alpha$ is played by the geodesic
representative of $\alpha$ in $N$.

\begin{proof}{}
Suppose the theorem is false; then there is a sequence of $\ep$-thick
geodesics $L_i$ parametrized by $[a_i,b_i]$, with foliations
$\alpha_i$, $\beta_i$ such that $|s_{\alpha_i}-s_{\beta_i}|\to\infty$
while $I_{s_\alpha}(\alpha,\beta)\to 0$.

Adjusting if necessary by automorphisms of $S$, we may assume that the
points $L_i(s_{\alpha_i})$ lie in a fixed compact subset of a
fundamental domain of the action of the mapping class group on
$\TT(S)$. Further, the quantity $I_t$ is invariant under scaling of
the measures, and so is defined on the compact space
$\PMF(S)\times\PMF(S)$. Thus we may take a convergent subsequence and
obtain a limit example $(L,\alpha,\beta)$ in which
$I_{s_\alpha}(\alpha,\beta) = 0$, and $|s_\beta-s_\alpha|=\infty$.
Without loss of generality we assume $s_\beta = +\infty$, so that $L$
must contain an $\ep$-thick ray $[s_\alpha,\infty)$. 

However, by lemma \ref{Thick implies complete}, this implies that the
foliation $\Phi_h$ is complete. Since $s_\beta = +\infty$, we have
$i(\beta,\Phi_h) = 0$, so that $\beta$ must be topologically
equivalent to $\Phi_h$. On the other hand, the fact that $s_\alpha \ne
+\infty$ implies that $i(\alpha,\Phi_h)>0$, and therefore that
$i(\alpha,\beta)>0$, a contradiction.
\end{proof}

Fix $\sigma\in\TT(S)$, and for $\alpha\in\MF(S)$
and $t\in \R$ define 
$$r_t(\alpha) = {e_t(\alpha)\over E_\sigma(\alpha)}.$$
This function is invariant under scaling, so it may be defined
for $[\alpha]\in\PMF(S)$. Since intersection number and extremal length
are continuous functions, it follows that the values
$r_{s_\alpha}(\alpha)$ depend continously on $[\alpha]$. Thus a
maximum is realized, and we can make these definitions:
\begin{eqnarray*}
\LL(\sigma,L) &=& \{[\lambda]\in\PMF(S): r_{s_\lambda}(\lambda) =
\max_{\alpha\in\MF(S)} r_{s_\alpha}(\alpha)\},\\
\SS(\sigma,L) &=& \{s_\lambda: [\lambda]\in\LL(\sigma,L)\}.
\end{eqnarray*}

The ratios $R_t(\alpha)$ are approximated by $r_t(\alpha)$, by virtue
of lemma \lref{Extremal length estimate}.

\medskip

We claim now that, for $[\lambda]\in\LL(\sigma,L)$, the ratio
$R_{s_\lambda}(\lambda)$ (or $r_{s_\lambda}(\lambda)$)
provides a good estimate for
$d(\sigma,L(s_\lambda))$. More precisely, 

\begin{lemma}{Ratio estimates distance}{}
There exists a constant $c_3(\ep,\chi(S))$ such that
$$R_{s_\lambda}(\lambda) \le \exp 2d(\sigma,L(s_\lambda) \le c_3R_{s_\lambda}(\lambda)$$
for $[\lambda]\in \LL(\sigma,L)$.
\end{lemma}

\begin{proof}{}
The left side of the inequality is immediate from (\ref{K ratio}).

Let
$[\mu]\in\PMF(S)$ denote the unique projective class of measured
foliations, such that  
$$d(\sigma,L(s_\lambda)) = \half\log R_{s_\lambda}(\mu),$$
and, in particular, $[\mu]$ maximizes $R_{s_\lambda}$ over $\PMF(S)$.
Let $A = R_{s_\lambda}(\mu) / R_{s_\lambda}(\lambda)$.
By virtue of (\ref{e is exp}), we have
$$r_{s_\lambda}(\mu) \le 2 e^{|s_\lambda-s_\mu|}r_{s_\mu}(\mu).$$
Since $R_{s_\lambda}(\mu) \le c_0r_{s_\lambda}(\mu)$, and by choice of
$\lambda$, $r_{s_\mu}(\mu)\le r_{s_\lambda}(\lambda)\le
R_{s_\lambda}(\lambda)$, we conclude 
$$A\le {c_0\over 2} e^{|s_\lambda-s_\mu|},$$
or
$$|s_\lambda-s_\mu| \ge \log A/2c_0.$$
Now suppose $A\ge 2c_0\exp D$, with $D$ the constant from lemma
\lref{Distance implies intersection}. It follows that
$$i(\lambda,\mu)^2\ge c_1E_{s_\lambda}(\lambda)E_{s_\lambda}(\mu).$$
On the other hand, lemma \ref{Extremal length and intersection number}
assures us that 
$$E_\sigma(\lambda)E_\sigma(\mu)\ge i(\lambda,\mu)^2.$$
Therefore we have
\begin{equation}
R_{s_\lambda}(\lambda)R_{s_\lambda}(\mu) \le {1\over c_1}
\label{R product small}
\end{equation}
or
$$A\le (c_1R^2_{s_\lambda}(\lambda))^{-1}.$$
To bound $A$ we need only note that there is a lower bound for
$R_{s_\lambda}(\lambda)$: There is a constant $\ell_0$ depending on
$\chi(S)$ such there is always a simple closed curve $\alpha$ on $S$
with $E_\sigma(\alpha)\le \ell_0$. Since $L$ is $\ep$-thick,
$R_t(\alpha)\ge \ep/\ell_0$ for all $t\in[a,b]$. 
Tracing through our definitions, 
it follows that $r_{s_\alpha}(\alpha) \ge \ep/\ell_0c_0$, and hence
that 
$R_{s_\lambda}(\lambda)\ge \ep/\ell_0c_0$ as well. This gives the desired
bound on $A$.
\end{proof}

\medskip

Next we will bound the diameter of
$\SS(\sigma,L)$. Suppose
$\lambda,\mu\in\LL(\sigma,L)$. If $|s_\lambda-s_\mu|\ge D$ then, as
above, inequality (\ref{R product small}) holds. On the other hand,
$R_{s_\lambda}(\lambda) \ge r_{s_\lambda}(\lambda)$, and 
$R_{s_\lambda}(\mu) \ge \half r_{s_\mu}(\mu)\exp|s_\lambda-s_\mu|$.
since, as in the previous paragraph,
$r_{s_\lambda}(\lambda)=r_{s_\mu}(\mu)\ge r_0$ for a fixed $r_0$,  we
get
$$\half r_0^2 e^{|s_\lambda-s_\mu|} \le {1\over c_1}$$
which bounds $|s_\lambda-s_\mu|$ by some $c_4$. 

Similarly, if $\lambda\in\LL(\sigma,L)$ and $\mu\in A_2(\sigma,L)$
then $inf_t R_t(\mu)$ is maximal in $\MF(S)$ and in paticular is at
least $\ep/\ell_0$. Thus we may repeat the argument of the previous
paragraph to obtain a bound on $|s_\lambda-s_\mu|$ in this case as
well.

Using (\ref{bound malpha}), we conclude that $T_2(\sigma,L)$ lies in a
bounded neighborhood of $\SS(\sigma,L)$, and thus its diameter is
bounded, and half the theorem is proven.

\medskip

To complete the proof, we will show that every point $t\in
T_1(\sigma,L)$ is at a bounded distance from $T_2(\sigma,L)$, or what is
now equivalent, from $\SS(\sigma,L)$.

Fix $t\in T_1(\sigma,L)$, and let
$[\mu]\in\PML(S)$ maximize $R_t$ -- so that
$$d(\sigma,L) = d(\sigma, L(t)) = \half \log R_t(\mu).$$
Since $d(\sigma,L(s_\lambda))\ge d(\sigma,L)$, and using
lemma \ref{Ratio estimates distance}, we obtain
$$c_3R_{s_\lambda}(\lambda) \ge R_t(\mu)$$
for any $\lambda\in \LL(\sigma,L)$. But by choice of $\mu$, $R_t(\mu)\ge
R_t(\lambda)$. Finally, $R_t(\lambda)\ge r_t(\lambda)\ge
\half\exp|t-s_\lambda|r_{s_\lambda}(\lambda)$, and
$R_{s_\lambda}(\lambda) \le c_0r_{s_\lambda}(\lambda)$. We conclude
$$c_0c_3 \ge \half \exp|t-s_\lambda|.$$
This bounds $|t-s_\lambda|$ by some $c_5$, and the proof is complete. 
\end{proof} % characterize projection

Theorem \ref{Quasi-map} follows too, recalling that the parameter $t$
measures half of Teichm\"uller arclength along $L$.
\end{proof} % quasi-map

We proceed to the contraction theorem:

\begin{theorem}{Contraction at a distance}{}
There exist constants $b_1,b_2$ depending on $\ep,\chi(S)$ such
that if $L$ is an $\ep$-thick geodesic and
$d(\sigma,L)> b_1$ then, setting $r=d(\sigma,L)-b_1$,
$${\rm diam}(\pi_L(\nbd_r(\sigma)))\le b_2.$$
\end{theorem}

\begin{proof}{}
Let $\tau$ be such that $d(\sigma,\tau)\le r.$ Suppose $\lambda\in\LL
(\sigma,L)$ and $\mu\in\LL(\tau,L)$.  By the proof of theorem
\ref{Characterize projection}, it is sufficient to find a
value for $b_1$ that allows us to bound $|s_\lambda-s_\mu|$.

Suppose $|s_\lambda - s_\mu|>D$, where $D$ is the constant in lemma
\ref{Distance implies intersection}. 
Then inequality (\ref{R product small}) holds again, or in other words
$$R_{\sigma,s_\lambda}(\lambda)R_{\sigma,s_\lambda}(\mu) \le {1\over c_1}.$$
By choice of $\lambda$, and lemma \ref{Ratio estimates distance},
$$R_{\sigma,s_\lambda}(\lambda) \ge c_6 \exp 
2d(\sigma,L),$$ 
for an appropriate constant $c_6$, and similarly with $\sigma$ and
$\lambda$ replaced by $\tau$ and $\mu$.
We note also that 
$$E_\sigma(\mu)\le e^{2r} E_\tau(\mu)$$
by choice of $\tau$ and by (\ref{K ratio}), and therefore
$$R_{\tau,s_\lambda}(\mu) \le e^{2r} R_{\sigma,s_\lambda}(\mu).$$
By definition of
$s_\mu$ and by lemma 
\ref{Extremal length estimate}, 
$$R_{\tau,s_\mu}(\mu)\le c_0 R_{\tau,s_\lambda}(\mu).$$
Putting all this together, we obtain
$$d(\sigma,L) + d(\tau,L) \le r + C$$
for some new constant $C$. 
By the definition of $r$, we obtain $d(\tau,L) \le C - b_1$. But, by
the choice of $\tau$ and the triangle inequality, 
$d(\tau,L)\ge b_1$, so if we choose $b_1 > C/2$ we obtain a
contradiction here, implying that
$|s_\lambda - s_\mu| \le D$. 
\end{proof}

\section{Applications}
\label{applications}

We mention some easily stated consequences of the projection theorems.
Some more substantive applications to the area of hyperbolic
three-manifolds will hopefully be presented in a future paper.

Define a {\em $K$-quasi-geodesic} in $\TT(S)$ as a path $\Gamma:I\to\TT(S)$,
say parametrized by arclength, such that for any $a,b\in I$, 
$$|a-b|\le K d(\Gamma(a),\Gamma(b)).$$

\begin{theorem}{Stability of quasi-geodesics}{}
Let $\Gamma$ be a $K$-quasi-geodesic path in $\TT(S)$, whose endpoints
are connected by an $\ep$-thick Teichm\"uller geodesic $L$. Then
$\Gamma$ remains in a $B(K,\ep)$-neighborhood of $L$. 
\end{theorem}

The proof is standard -- see \cite{minsky-slow} for a proof in a very
similar situation. The idea is that any segment of $\Gamma$ which is
outside a sufficiently large neighborhood of $L$ is, because of 
the contractive properties of $\pi_L$, much less efficient than 
the path obtained by going back to $L$ and moving along the
projection. 

\medskip

Fix a non-peripheral simple closed curve $\alpha$ in $S$. The set 
$$Thin(\alpha,\delta) = \{\sigma\in\TT(S): E_\sigma(\alpha)\le\delta\}$$
is somewhat analogous to a horoball in hyperbolic space.
We have the following:

\begin{theorem}{Cusp projections}{}
Given $\ep, \delta$ there is a constant $B$ such that 
$${\rm diam}(\pi_L(Thin(\alpha,\delta)))\le B$$
for any simple closed curve $\alpha\subset S$ and $\ep$-thick geodesic
$L$. 
\end{theorem}

\begin{proof}{}
We show first that there is some $\delta_0$ and $B$ for which the bound
${\rm diam}(\pi_L(Thin(\alpha,\delta_0)))\le B$ holds. In fact we will show
that $\pi_L(Thin(\alpha,\delta_0))$ lies in a bounded neighborhood of
$L(s_\alpha)$. 

By the results of the previous section, it suffices to bound
$|s_\lambda-s_\alpha|$ for $\lambda\in\LL(\sigma,L)$, when
$\sigma\in Thin(\alpha,\delta_0)$. Suppose $|s_\lambda-s_\alpha|\ge
D$. Then, as in (\ref{R product small}) we get
$$R_{s_\alpha}(\alpha)R_{s_\alpha}(\lambda) \le {1\over c_1}.$$
Now, $R_{s_\alpha}(\alpha) \ge \ep/\delta_0$, so
$R_{s_\alpha}(\lambda) \le \delta_0/\ep c_1$.
On the other hand, $R_{s_\alpha}(\lambda)\ge r_{s_\lambda}(\lambda)
\ge r_0$ as in the previous section, so if we choose $\delta_0 < \ep
c_1 r_0$ we obtain a contradiction, and conclude $|s_\lambda -
s_\alpha|< D$.

\medskip

To prove the theorem for general $\delta > \delta_0$, we observe that
$Thin(\alpha,\delta)$ lies in a neighborhood of
$Thin(\alpha,\delta_0)$ of radius $\half\log \delta/\delta_0$ (since
the extremal length of $\alpha$ can always be reduced at a constant
rate along an appropriate Teichm\"uller geodesic), and apply theorem
\ref{Contraction at a distance}. 

Let $\tau$ be a point at a distance $a$ from
$Thin(\alpha,\delta_0)$, and let $\sigma\in Thin(\alpha,\delta_0)$ be
such that $d(\sigma,\tau) = a$. We may assume $a>b_1$. If $d(\sigma,L)
> a + b_1$ then by theorem \ref{Contraction at a distance}
we have $d(\pi_L(\sigma),\pi_L(\tau)) \le b_2$ (here distance between
sets is the maximal distance between pairs of elements), and we are done.

If $d(\sigma,L) \le a + b_1$ then $d(\tau,\pi_L(\sigma)) \le 2a +
b_1$. This implies $d(\tau,\pi_L(\tau)) \le 2a + b_1$, so
$d(\pi_L(\sigma),\pi_L(\tau))\le 4a + 2b_1$, and again we are done.

\end{proof}


\ifx\undefined\bysame
\newcommand{\bysame}{\leavevmode\hbox to3em{\hrulefill}\,}
\fi
\begin{thebibliography}{10}

\bibitem{abikoff}
W.~Abikoff, {\em The real-analytic theory of {T}eichm\"uller space},
  Springer-Verlag, 1980, Lecture Notes in Mathematics no. 820.

\bibitem{bers:pseudoanosov}
L.~Bers, {\em An extremal problem for quasiconformal mappings and a theorem by
  {T}hurston}, Acta Math. {\bf 141} (1978), 73--98.

\bibitem{travaux}
A.~Fathi, F.~Laudenbach, and V.~Poenaru, {\em Travaux de {T}hurston sur les
  surfaces}, vol. 66-67, Asterisque, 1979.

\bibitem{gardiner}
F.~Gardiner, {\em {T}eichm\"{u}ller theory and quadratic differentials}, Wiley
  Interscience, 1987.

\bibitem{gardiner-masur}
F.~Gardiner and H.~Masur, {\em Extremal length geometry of {T}eichm\"uller
  space}, Complex Variables, Theory and Applications {\bf 16} (1991), 209--237.

\bibitem{hubbard-masur}
J.~Hubbard and H.~Masur, {\em Quadratic differentials and foliations}, Acta
  Math. {\bf 142} (1979), 221--274.

\bibitem{kerckhoff}
S.~Kerckhoff, {\em The asymptotic geometry of {T}eichm\"uller space}, Topology
  {\bf 19} (1980), 23--41.

\bibitem{kerckhoff:nielsen}
\bysame, {\em The {N}ielsen realization problem}, Ann. of Math. {\bf 117}
  (1983), 235--265.

\bibitem{levitt}
G.~Levitt, {\em Foliations and laminations on hyperbolic surfaces}, Topology
  {\bf 22} (1983), 119--135.

\bibitem{masur:teichgeo}
H.~A. Masur, {\em On a class of geodesics in {T}eichm\"uller space}, Ann. of
  Math. {\bf 102} (1975), 205--221.

\bibitem{masur:uergodic}
\bysame, {\em Uniquely ergodic quadratic differentials}, Comment. Math. Helv.
  {\bf 55} (1980), 255--266.

\bibitem{masur:geoflow}
\bysame, {\em Transitivity properties of the horocyclic and geodesic flows on
  moduli space}, J. D'analyse Math. {\bf 39} (1981), 1--10.

\bibitem{minsky:extremal}
Y.~Minsky, {\em Extremal length estimates and product regions in
  {T}eichm\"uller space}, Stony Brook IMS Preprint \#1994/11.

\bibitem{minsky:2d}
\bysame, {\em Harmonic maps, length and energy in {T}eichm\"uller space}, J. of
  Diff. Geom. {\bf 35} (1992), 151--217.

\bibitem{minsky:slowmaps}
\bysame, {\em Teichm\"uller geodesics and ends of hyperbolic 3-manifolds},
  Topology {\bf 32} (1993), 625--647.

\bibitem{wolpert:noncomplete}
S.~A. Wolpert, {\em Noncompleteness of the {Weil-Peterson} metric for
  {Teichm\"uller} space}, Pacific J. of Math. {\bf 61} (1975), 573--577.

\bibitem{wolpert:nielsen}
\bysame, {\em Geodesic length functions and the {N}ielsen problem}, J.
  Differential Geom. {\bf 25} (1987), 275--296.

\end{thebibliography}


\begin{thebibliography}{BGS}

\bibitem[A]{abikoff}
W. Abikoff.
\newblock {\it The Real-Analytic Theory of Teichm\"uller Space}.
\newblock {\it Lecture Notes in Mathematics no. 820}, Springer-Verlag, 1980.

\bibitem[B]{bers}
L. Bers.
\newblock An extremal problem for quasiconformal mappings and a
	  theorem by Thurston,
\newblock {\it Acta Math.} {\bf 141} (1978), 73--98.

\bibitem[FLP]{travaux}
A. Fathi, F. Laudenbach, and V. Poenaru.
\newblock {\it Travaux de Thurston sur les Surfaces}.
\newblock Volume~{\bf 66-67}, Asterisque, 1979.

\bibitem[Ga]{gardiner}
F. P. Gardiner. 
\newblock {\it Teichm\"uller Theory and Quadratic Differentials.} 
\newblock Wiley-Interscience, 1987.

\bibitem[GM]{gardiner-masur}
F. Gardiner and H. Masur.
\newblock Extremal length geometry of Teichm\"uller space,
\newblock {\it preprint.}

\bibitem[HM]{hubbard-masur}
J. Hubbard and H. Masur.
\newblock Quadratic differentials and foliations,
\newblock {\it Acta Math.} {\bf 142} (1979), 221--274


\bibitem[K1]{kerckhoff}
S. Kerckhoff.
\newblock {The asymptotic geometry of Teichm\"uller space}.
\newblock {\it Topology} {\bf 19} (1980), 23--41.

\bibitem[K2]{kerckhoff2}
S. Kerckhoff.
\newblock {The Nielsen realization problem}.
\newblock {\it Annals of Math.} {\bf 117} (1983), 235--265.

\bibitem[L]{levitt}
G. Levitt.
\newblock {Foliations and laminations on hyperbolic surfaces}.
\newblock {\it Topology} {\bf 22} (1983), 119--135.

\bibitem[Ma]{masur}
H. Masur.
\newblock On a class of geodesics in Teichm\"uller space,
\newblock {\it Annals of Math.} {\bf 102} (1975), 205--221.

\bibitem[Mi1]{minsky-2d}
Y. Minsky. 
\newblock Harmonic maps, length and energy in Teichm\"uller space,
\newblock To appear in {\it J. of Diff. Geom.}

\bibitem[Mi2]{minsky-slow}
Y. Minsky.
\newblock Teichm\"uller geodesics and ends of hyperbolic 3-manifolds,
\newblock {\it To appear.}

\bibitem[W]{wolpert}
S. Wolpert.
\newblock Geodesic length functions and the Nielsen problem.
\newblock {\it J. of Diff. Geom.} {\bf 25} (1987), 275--296.


\end{thebibliography}
\end{document}